\documentclass[oneside,english]{amsart}
\usepackage[T1]{fontenc}
\usepackage[latin9]{inputenc}
\usepackage{amsthm}

\makeatletter
\numberwithin{equation}{section}
\numberwithin{figure}{section}
\theoremstyle{plain}

\theoremstyle{definition}

\theoremstyle{plain}

\theoremstyle{plain}

\theoremstyle{remark}

\makeatother

\usepackage{babel}
\providecommand{\definitionname}{Definition}
\providecommand{\lemmaname}{Lemma}
\providecommand{\propositionname}{Proposition}
\providecommand{\remarkname}{Remark}
\providecommand{\theoremname}{Theorem}

\begin{document}
\title{Quasiconformal Whitney Partition}
\author{Vladimir Gol'dshtein$^{1}$ and Nahum Zobin$^{2}$} 
 \footnote { Ben Gurion University of the Negev} 
 \footnote { College of William and Mary, VA, USA}

\begin{abstract}
Whitney partition is a very important concept in modern analysis. We
discuss here a quasiconformal version of the  Whitney partition that can
be usefull for Sobolev spaces. 
\end{abstract}

\maketitle

\section{Introduction}
 We start with a definition of a Whitney partition of domains in $\mathrm{\mathbf{R}^{n}}$.
Classical Whitney partition is a partition of a bounded domain $\Omega$
into diadic cubes with disjoint interiors and edges comparable to the distance
to $\partial\Omega$. Its modern generalization that is called a Whitney
partition is a partition into convex polyhedra with similar properties.
Let us give a more accurate definition (see, for example \cite{R}).

Let $\Omega$ be a bounded domain of $\mathrm{\mathbf{R}^{n}}$ and
$\Lambda=\mathrm{\mathbf{R}^{n}}\setminus\Omega$. Let $E_{i}$ be
a family of convex closed $n$-dimensional polyhedra in $\mathrm{\mathbf{R}^{n}}$,
disjoint from $\Lambda$, covering $\Omega$ and with pairwise disjoint
interiors. We will also need these polyhedra to have uniformly bounded
ratio $K_{A}$ of their exterior to interior radii. By the interior
radius of a set $A$ with nonempty interior we mean the greatest radius
$r_{A}$ of a ball contained in $A$; similarly the exterior radius of
a set $A$ is the smallest radius $R_{A}$ of a ball containing $A$.
The ratio $K_{A}$ will be called the dilatation of $A$.

Such families of polyhedra with uniformly bounded dilatation will
be called uniformly regular. We will also demand that the edges of
these polyhedra are long, i.e. they have lengths uniformly comparable to
the diameter of the polyhedron. As an example we can take a family
of  dyadic cubes disjoint from $\Lambda$. The diameter of
$E_{i}$ will be denoted by $\delta\left(E_{i}\right)$. The set $\Lambda$
will be called the residual set of the family $\left\{ E_{i}\right\} $
and $K_{int}:=\sup_{i}K_{E_{i}}$ will be called the interior dilatation
of the family $\left\{ E_{i}\right\} $.

There are two conditions which are usually imposed on such families:
\begin{enumerate}
\item $\delta\left(E_{i}\right)\geq C\mathrm{\:dist}\left(E_{i},\Lambda\right)$,
\item $\delta\left(E_{i}\right)\leq C\mathrm{^{-1}\:dist}\left(E_{i},\Lambda\right)$. 
\end{enumerate}
If a family $E_{i}$ satisfies both of these conditions, it is called a
Whitney family.

We generalize this definition for a more general setting under additional assumption that the residual set
$\mathrm{\mathbf{R}^{n}}\setminus\Omega'$ is unbounded and connected. We will call such domains as simple domains. 

{\bf Definition.}(Rough Whitney family). {\it Let $\Omega$ be a bounded subset of $\mathrm{\mathbf{R}^{n}}$,
$\Lambda=\mathrm{\mathbf{R}^{n}}\setminus\Omega$ and $\left\{ E_{i}\right\} $
be a family of closed sets in $\mathrm{\mathbf{R}^{n}}$ with nonempty
interiors $\mathrm{V_{i}=Int}(E_{i})$ and the closures of $V_{i}$
coinciding with $E_{i}$. The sets $E_{i}$ are disjoint from $\Lambda$,
covering $\Omega$ and their interiors $V_{i}$ are pairwise disjoint.
The family $\left\{ E_{i}\right\} $ has bounded interior dilatation.

The family $E_{i}$ should also satisfy  two geometric regularity conditions
for some constant $C>0$:

1.$\delta\left(E_{i}\right)\geq C\:\mathrm{dist}\left(E_{i},\Lambda\right)$,

2.$\delta\left(E_{i}\right)\leq C^{-1}\:\mathrm{dist}\left(E_{i},\Lambda\right)$. }

If a family $E_{i}$ satisfies all these conditions, we will call it as
a rough Whitney family.

The smallest possible constant $C$ will be called the exterior dilatation
$K_{ext}$ of the rough Whitney family.

Our main result connects concepts of Whitney families and rough Whitney
families.

{\bf Theorem.}
{\it For any simple domain $\Omega$ quasiconformal image of any its Whitney family is a rough Whitney family.}

More detailed formulation with corresponding estimates will be discussed later. 

We prove this result with the help of classical estimates of a corresponding capacity.

A two dimensional version of this result can be found in \cite{KRZ}.

\section{Quasiconformal Whitney family }

A well-ordered triple $(F_{0},F_{1};\Omega)$ of nonempty sets, where
$\Omega$ is an open set in $\mathrm{\mathbf{R^{n}}}$, and $F_{0}$,
$F_{1}$ are closed subsets of $\overline{\Omega}$, is called a condenser
on the Euclidean space $\mathrm{\mathbf{R^{n}}}$.

The value 
\[
\mathrm{Cap_{p}}(E)=\mathrm{Cap_{p}}(F_{0},F_{1};\Omega)=\inf\int\limits _{\Omega}|\nabla v|^{p}dx,
\]
where the infimum is taken over all Lipschitz non-negative functions
$v:\overline{\Omega}\to\mathrm{\mathbf{R}}$, such that $v=0$ on
$F_{0}$, and $v=1$ on $F_{1}$, is called $p$-capacity of the condenser
$E=(F_{0},F_{1};\Omega)$. For $p=n$ it is the classical conformal
capacity. We will use notation $\mathrm{Cap}(F_{0},F_{1};\Omega)$
for the conformal capacity. Of course the set of admissible functions
can be empty. In this case $\mathrm{Cap_{p}}(F_{0},F_{1};\Omega)=\infty$.

For $1<p$ for  a finite value of $p$-capacity $0\leq\mathrm{Cap}_{p}(F_{0},F_{1};\Omega)<+\infty$
there exists a unique continuous weakly differentiable function $u_{0}$
(an extremal function) such that: 
\[
\mathrm{Cap}_{p}(F_{0},F_{1};\Omega)=\int\limits _{\Omega}|\nabla u_{0}|^{p}dx.
\]

{\bf Definition.}
{ \it Let $\varphi:\Omega\rightarrow\Omega\prime$ be a homeomorphism between
two domains in $\mathrm{\mathbf{R^{n}}}$, $n\ge2$. Then $\varphi$
is said to be $Q$-quasiconformal, $Q\ge1$, if 

\[
Q^{-1}\mathrm{Cap} (F_{0},F_{1};\Omega)\leq\mathrm{Cap}(\varphi\left(F_{0}\right),\varphi\left(F_{1}\right);\Omega')\leq Q\:\mathrm{Cap}(F_{0},F_{1};\Omega)
\]
for any condenser $E=(F_{0},F_{1};\Omega)$.

The minimal possible constant $Q$ will be called the (quasiconformal)
dilatation of $\varphi$.}

This geometric definition is a global requirement that quickly yields
many important properties of quasiconformal mappings. For example,
the inverse of a quasiconformal mapping is automatically quasiconformal,
quasiconformal mappings are weakly differentiable and its weak derivatives are
integrable in degree $n$, etc...

For any $Q$-quasiconfromal homeomorphism $\varphi:\mathrm{\mathbf{R^{n}}}\to\mathrm{\mathbf{R^{n}}}$
images of closed balls have uniformly bounded interior dilatation
$K_{int}$ that depends only on $n$ and $Q$ \cite{Vais}. Any $Q$-quasiconformal
image of closed balls $\bar{B}(0,r)$ will be called the closed $Q$-quasiconformal
ball or $Q$-quasiball. For example any convex polyhedra is a closed
$Q$-quasiconfromal ball. The constant $Q$ depends on $n$ and its
interior dilatation. This collection of known facts can be formalized
as 

{\bf Proposition 1.} (\cite{Vais}) {\it For any $Q$-quasiconfromal homeomorphism $\varphi:\mathrm{\mathbf{R^{n}}}\to\mathrm{\mathbf{R^{n}}}$
and any closed ball $B:=B(x_{0},r)$ the interior dilatation $K_{B,\Omega}$
is less or equal to a constant $C(Q,n)$ that depends only on the
dilatation $Q$ and $n$. }

This proposition will be generalized for balls in domains under an
additional conditions for balls.

We will use notations: $B(x_{0},r)$ is an open ball of radius $r$
and center $x_{0}$ and $\bar{B}(x_{0},r)$ is its closure.
 
{\bf Definition.}
{\it Let $\Omega$ be a domain in $\mathrm{\mathbf{R^{n}}}$ and a closed
ball $\overline{B}(x_{0},r)\subset\Omega$. Choose such concentric
open ball $B(x_{0},K_{B(x_{0},r),\Omega}r)\subset\Omega$ that its
closure $\bar{B}(x_{0},K_{B(x_{0},r),\Omega}r)$ intersects $\partial\Omega$.
We will call the constant $K_{B(x_{0},r),\Omega}\in(1,\infty]$ the
embedding coefficient of the ball $\bar{B}(x_{0},r)$ in $\Omega$.}

The embedding coefficient of any ball in $\mathbf{R^{n}}$ is $\infty$.
In section 4 we will prove the following property of quasiconformal
homeomorphisms.

{\bf Proposition 2.}
{\it Let $\varphi:\Omega\to\mathrm{\mathbf{R}^{n}}$ be a $Q$-quasiconformal
homeomorphism of a domain $\Omega\subset\mathrm{\mathbf{R}^{n}}$
and \textup{$\overline{B}(x_{0},r)$} is a closed ball with the embedding
coefficient $K_{B(x_{0},r),\Omega}$. Then the interior dilatation
$K_{\varphi(B(x_{0},r))}$ of its image $\varphi(B(x_{0},r)$ depends
only on $Q,K_{B(x_{0},r),\Omega}$ and $n$.}

We will call a closed set $E$ the closed relative $Q$-quasiball
if there exist a domain $\Omega\subset\mathrm{\mathbf{R}^{n}}$, a
closed ball $\overline{B}(x_{0},r)$ and a $Q$-quasiconformal homeomorphism
$\varphi:\Omega\to\mathrm{\mathbf{R}^{n}}$ such that $E=\varphi\left(\overline{B}(x_{0},r)\right)$.
We let $K_{E}$ denote the embedding coefficient $K_{B(x_{0},r),\Omega}$.

{\bf Definition} ($Q$-quasiconfromal Whitney family). {\it Let $\Omega$ be a bounded domain
in $\mathrm{\mathbf{R}^{n}}$ and $E_{i}$ be a family of closed relative
$Q$-quasiballs in $\mathrm{\mathbf{R}^{n}}$. The sets $E_{i}$ are
disjoint from $\Lambda$, covering $\Omega$, with pairwise disjoint
interiors $V_{i}$ and their embedding coefficients $K_{i}:=K_{E_{u}}$ are
uniformly bounded.

If a family $E_{i}$ satisfies this conditions, it is called a $Q$-quasiconformal
Whitney family.}

Because any convex polyhedra is a quasiconformal image of the unit
ball, any Whitney family is a $Q$-quasiconformal Whitney family.
Of course a rough Whitney family is not necessary a $Q$-quasiconformal
Whitney family. 

\section{Classical estimates of conformal capacity}

We will need two estimates of the conformal capacity (see, for example
\cite{Vais,GR}). Because proofs are very simple we reproduce here
slight modifications of both estimates, using classical embedding
theorems. 

Choose an open ball $B(x_{0},r)$ an open ball of radius $r$ and center
$x_{0}$ and a positive constant $C_{0}>1$. Denote $R:=C_{0}\:r$.

{\bf Lemma 1.} 
{\it $\mathrm{Cap}\left(\bar{B}(x_{0},r),\mathrm{\mathbf{R}^{n}}\setminus B(x_{0},R);\mathrm{\mathbf{R}^{n}}\right)=\omega_{n-1}\left[\ln\left(C_{0}\right)\right]^{1-n}.$ \\
Here $\omega_{n-1}$is the volume of the unit $\left(n-1\right)$-dimensional
sphere $S^{n-1}$.}

{\bf Proof.}
Let $u:\mathrm{\mathbf{R}^{n}}\to\mathrm{\mathbf{R}}$ be an admissible
smooth function for condensor $E=\left(\bar{B}(x_{0},r),\mathrm{\mathbf{R}^{n}}\setminus B(x_{0},R);\mathrm{\mathbf{R}^{n}}\right)$.
By definition of an admissible function $u(x)\equiv0$ on $\bar{B}(x_{0},r)$
and $u(x)\equiv1$ on $\mathrm{\mathbf{R}^{n}}\setminus B(x_{0},R)$.
Using H$\ddot{\mathrm{o}}$lder inequality for the spherical coordinate
system $(\rho,\sigma)$ we have
\[
1\leq\int_{r}^{R}\nabla u(\rho,\sigma)d\rho\leq\left(\int_{r}^{R}\left|\nabla u(\rho,\sigma)\right|^{n}\rho^{n-1}d\rho\right)^{1/n}\left(\int_{r}^{R}\rho^{-1}d\rho\right)^{\frac{n-1}{n}}.
\]

Here $\rho:=|x|$ and $\sigma$ are spherical coordinates on the unit
sphere $S^{n-1}:=S^{n-1}(0,1)$. By previos calculations

\[
\frac{1}{\left[\ln\left(C_{0}\right)\right]^{n-1}}=\frac{1}{\left[\ln\left(\frac{R}{r}\right)\right]^{n-1}}\leq\int_{r}^{R}\left|\nabla u(\rho,\sigma)\right|^{n}\rho^{n-1}d\rho
\]
for any $\sigma\in S^{n-1}$.

Integrating this inequality over $S^{n-1}$ we obtain finally
\[
\frac{\omega_{n-1}}{\left[\ln\left(C_{0}\right)\right]^{n-1}}\leq\int_{S^{n-1}}\int_{r}^{R}\left|\nabla u(\rho,\sigma)\right|^{n}\rho^{n-1}d\rho=\int_{\mathrm{\mathbf{R}}^{n}}\left|\nabla u(x)\right|^{n}dx
\]
where $\omega_{n-1}$ is area of $S^{n-1}$.

The function $u_{0}(x):=\ln\left(\frac{R}{\left|x\right|}\right)\left[\ln\left(C_{0}\right)\right]^{-1}$
for $r\leq\left|x\right|\leq R$ , $0$ for any $\left|x\right|\geq R$
and $1$ for any $\left|x\right|\leq r$ is the extremal function
for condensor $E$. It follows from direct calculations: 
\[
\int_{\mathrm{\mathbf{R}}^{n}}\left|\nabla u_{0}(x)\right|^{n}dx=\omega_{n-1}\left[\ln\left(C_{0}\right)\right]^{-n}\int_{r}^{R}\rho^{-1}d\rho=\omega_{n-1}\left[\ln\left(C_{0}\right)\right]^{1-n}.
\]

Therefore
\[
\mathrm{Cap}\left(\bar{B}(x_{0},r),\mathrm{\mathbf{R}^{n}}\setminus B(x_{0},R);\mathrm{\mathbf{R}^{n}}\right)=\omega_{n-1}\left[\ln\left(C_{0}\right)\right]^{1-n}.
\]

We will call a connected closed set a continuum. Choose two concentric
$(n-1)$-dimensional spheres $S^{n-1}(0,r)$ and $S^{n-1}(0,R)$,
$R\geq r$ and two continua $F_{0}$ and $F_{1}$ that join spheres.
We shall use notations $R=C_{0}\:r$, $C_{0}>1$ and $D_{r,R}:=\bar{B(0,R)\setminus B(0,r)}$.
In the next lemma we will prove a below estimate of conformal capacity
of condensor $E=(F_{0},F_{1},\mathrm{\mathbf{R}^{n})}$ using embedding
theorems for the unit sphere. 

{\bf Lemma 2.}  {\it$\mathrm{Cap}\left(F_{0},F_{1};\mathrm{\mathbf{R}^{n}}\right)\geq C(n)\:\ln\left(C_{0}\right)$,
where a constant $C(n)$ depends on $n$ only.}

{\bf Proof.}
Because continua $F_{0},F_{1}$ join spheres $S^{n-1}(0,r)$ and
$S^{n-1}(0,R)$ intersections $F_{0,\rho}:=S_{\rho}^{n-1}\cap F_{0}$
and $F_{1,\rho}:=S_{\rho}^{n-1}\cap F_{1}$ are not empty for any
sphere $S_{\rho}^{n-1}:=S^{n-1}(0,\rho)$, $r\leq\rho\leq R$. Any  function
$u$
admissible for the conformal capacity of the condensor $E$
 is also admissible for conformal capacity of any condensor $E_{\rho}=(F_{0,\rho},F_{1,\rho},\mathrm{S_{\rho}^{n-1})}$
in the sphere $S_{\rho}^{n-1}$. By elementary calculations in the
spherical coordinates $\rho,\sigma$ we have 
\[
\int_{r}^{R}\left(\int_{S_{\rho}^{n-1}}\left|\nabla u(\sigma,\rho)\right|^{n}d\sigma\right)d\rho\leq\int_{D_{R,r}}\left|\nabla u(x)\right|^{n}dx\leq\int_{\mathrm{\mathbf{R}^{n}}}\left|\nabla u(x)\right|^{n}dx
\]
where $D_{R,r}$ is a closed ring between spheres $S^{n-1}(0,r)$
and $S^{n-1}(0,R)$.

Using similarities $\varphi_{\rho}(x)=\rho x$ we get 
\[
\int_{S_{\rho}^{n-1}}\left|\nabla u(\sigma,\rho)\right|^{n}d\sigma=\frac{1}{\rho}\int_{S^{n-1}(0,1)}\left|\nabla\tilde{u}(\sigma)\right|^{n}d\sigma,
\]
where $\tilde{u}(x)=u(\varphi_{\rho}(x))$.

By the classical Sobolev inequality for the unit sphere $S^{n-1}$
we have 
\[
1=\left\Vert \tilde{u}|L_{\infty}(S^{n-1})\right\Vert ^{n}\leq K(n)\int_{S^{n-1}(0,1)}\left|\nabla\widetilde{u}(\sigma)\right|^{n}d\sigma
\]
where constant $K(n)$ depends on $n$ only. 

Combining these two estimates we obtain 
\[
1=\left\Vert u|L_{\infty}(S_{\rho}^{n-1})\right\Vert ^{n}=\left\Vert \tilde{u}|L_{\infty}(S^{n-1}(0,1))\right\Vert ^{n}\leq\rho K(n)\int_{S_{\rho}^{n-1}}\left|\nabla u(\sigma,\rho)\right|^{n}d\sigma.
\]
Dividing by $\rho$ and integrating we finally get
\[
\frac{\ln\left(C_{0}\right)}{K(n)}=\frac{1}{K(n)}\int_{r}^{R}\frac{d\rho}{\rho}\leq\int_{r}^{R}\left(\int_{S_{\rho}^{n-1}}\left|\nabla u(\sigma,\rho)\right|^{n}d\sigma\right)d\rho
\]
\[
\leq\int_{\mathrm{\mathbf{R}^{n}}}\left|\nabla u(x)\right|^{n}dx
\]
for any admissible function $u$ of condensor $E$. By definition
of conformal capacity

\[
\frac{\ln\left(C_{0}\right)}{K(n)}\leq\mathrm{Cap}\left(F_{0},F_{1};\mathrm{\mathbf{R}^{n}}\right).
\]

\section{Local estimates of dilatations}

Let $\varphi:\Omega\to\Omega'$ be a $Q$-quasiconformal homeomorphism
of a domain $\Omega\subset\mathrm{\mathbf{R}^{n}}$ onto a domain
$\Omega'\subset\mathrm{\mathbf{R}^{n}}$. Choose a closed ball $\bar{B}(x_{0},r)\subset\Omega$
with the embedding coefficient $C_{r}:=K_{B\{x_{0},r),\Omega}$ .

Our goal is to prove for $F_{r}:=\varphi(\bar{B}(x_{0},r))$ the following
two inequalities:
\begin{enumerate}
\item There exists a constant $C\mathrm{_{1}(Q,C_{r},n)}$ such that 
$$K_{int}(F_{r})\leq C\mathrm{_{1}(Q,C_{r},n)}$$
\item There exists a constant $C\mathrm{_{2}(Q,C_{r},n)}$ such that 
$$  C^{-1}\:\mathrm{dist}\left(F_{r},\Omega'\right)\leq \delta\left(F_{r}\right)\leq C\:\mathrm{dist}\left(F_{r},\Omega'\right)$$,
\end{enumerate}
Constants $C_{1},C_{2}$ depend only on $Q,C_{r},n$ and do not
depend on choice of domains $\Omega,\Omega'$ and the $Q$-quasiconformal
homeomorphism $\varphi$. 

{\bf Proposition 3.}  {\it Suppose that $\Omega$,$\Omega'$ are domains
in $\mathbf{R}^{n}$ and the residual set $\Lambda=\mathbf{R}^{n}\setminus\Omega$
of $\Omega$ is unbounded and connected. Let $\varphi:\Omega\to\Omega'$
be a $Q$-quasiconformal homeomorphism of a domain $\Omega\subset\mathrm{\mathbf{R}^{n}}$
onto a domain $\Omega'\subset\mathrm{\mathbf{R}^{n}}$ and $\bar{B}(x_{0},r)\subset\Omega$
be a closed ball with the embedding coefficient $K_{B(x_{0},r),\Omega}$
. Then there exists such positive constant $C_{1}(Q,C_{r},n)$ that
$K_{int}(F_{r})\leq C_{1}(Q,K_{B(x_{0},r),\Omega},n)$.

The constant \textup{$C_{1}(Q,K_{B(x_{0},r),\Omega},n)$ depends on
$Q,K_{B(x_{0},r),\Omega}$ and $n$ only.}}

{\bf Proof.}
Denote by $\bar{B}_{r}:=\bar{B}(y_{0},\bar{r})\subset F_{r}$ a closed
ball of a maximal radius that belongs to $F_{r}:=\varphi(B(x_{o},r))$
and has the center at a point $y_{0};=\varphi(x_{0})$, by $B_{R}:=B(y_{0},\bar{R})$
an open ball of a minimal radius whose closure contains $F_{r}$,
by $CB_{R}:=\mathrm{\mathbf{R}^{n}}\setminus B_{R}$, by $K_{r}$
the interior dilatation $K_{int}(F_{r})$ and by $C_{r}$ the embedding
coefficient $K_{B(x_{0},r),\Omega}$. By definition of the interior
dilatation $K_{r}\leq\frac{\bar{R}}{\bar{r}}$. 

By definition of the conformal capacity and previous Lemma 
\begin{equation}
\mathrm{Cap}(B_{R},\bar{\Omega}'\setminus B,\Omega')\leq\mathrm{Cap}(B_{r},CB_{R},\mathrm{\mathbf{R}^{n}})=\omega_{n-1}\left[\ln\left(\frac{\bar{R}}{\bar{r}}\right)\right]^{1-n} 
\end{equation}
\begin{equation}
\leq\omega_{n-1}\left[\ln\left(K_{r}\right)\right]^{1-n}.\label{eq:zero}
\end{equation}

Because a homeomorphism $\varphi^{-1}$ is $Q$-quasiconformal we
have for two compact sets $F_{0}:=\varphi^{-1}(\bar{B}_{r})\subset\bar{B}(x_{0},r)$
and $F_{1}=\varphi^{-1}(\bar{\Omega}'\setminus B_{R})\subset\Omega\setminus B(x_{0},r)$
the following inequality 
\begin{equation}
\mathrm{Cap}(F_{0},F_{1},\Omega)\leq Q\mathrm{\:Cap}(\bar{B}_{r},\bar{\Omega}'\setminus B_{R},\Omega')\leq\omega_{n-1}\left[\ln\left(K_{r}\right)\right]^{1-n}.\label{eq:first}
\end{equation}

Let us estimate $\mathrm{Cap}(F_{0},F_{1},\Omega)$ with the help
of Lemma 2. By construction both sets $F_{0},F_{1}$
have nonempty intersections $S_{0},S_{1}$ with the sphere $S(x_{0},r)$.
We distinguish two different cases: 

$\mathrm{1.\;dist}(S_{0},S_{1})\leq\min\left(\left(C_{r}-1\right)\frac{r}{2},\frac{r}{2}\right)$;

$\mathrm{2.\;dist}(S_{0},S_{1})>\min\left(\left(C_{r}-1\right)\frac{r}{2},\frac{r}{2}\right)$.

Let us use a short notation $\bar{r}:=\min\left(\left(C_{r}-1\right)\frac{r}{2},\frac{r}{2}\right).$

Choose points $y_{o}\in S_{0}$ and $\bar{y_{1}\in S_{1}}$ such that
$\mathrm{dist}(S_{0},S_{1})=|y_{1}-y_{0}|$. Let $\tilde{y}:=y_{0}+\frac{y_{1}-y_{0}}{2}$
and $B_{1}:=\bar{B(}\tilde{y},\bar{r})$ be a closed ball with center
at $\tilde{y}$. 

In the first case $|y_{1}-y_{0}|\leq\bar{r}$ and continua $F_{0},F_{1}$
intersect any sphere $S(\tilde{y},\rho)$ for any $\bar{r}\leq\rho\leq2\bar{r}$
. By Lemma 2

\begin{equation}
\mathrm{\mathrm{Cap}\left(F_{0},F_{1};\Omega\right)\geq\mathrm{Cap}\left(F_{0}\cap B_{1},F_{1}\cap B_{1};B_{1}\right)}\geq C(n)\:\ln2.\label{eq:second}
\end{equation}

Combining this inequality and inequality \ref{eq:first} we obtain
finally
\begin{equation}
C(n)\:\ln2\leq\omega_{n-1}\left[\ln\left(K_{r}\right)\right]^{1-n}.\label{eq:final1}
\end{equation}

The proposition is proved for the first case when
$$\mathrm{\;dist}(S_{0},S_{1})\leq\min\left(\left(C_{r}-1\right)\frac{r}{2},\frac{r}{2}\right).$$

In the second case$\mathrm{\;dist}(S_{0},S_{1})\leq\min\left(\left(C_{r}-1\right)\frac{r}{2},\frac{r}{2}\right)$ \\
we have $|y_{1}-y_{0}|>\bar{r}$. By definition of capacity
\[
\mathrm{Cap}\left(F_{0},F_{1};\Omega\right)\geq\mathrm{Cap}\left(F_{0}\cap B(x_{0},C_{r}r),F_{1}\cap B(x_{0},C_{r}r);B(x_{0},C_{r}r)\right).
\]

There exists a quasiconformal homeomorphism $\psi$ that maps $B(x_{0},C_{r}r)$
onto itself, maps any sphere $S(x_{0},\rho)$, $0<\rho<C_{r}r$ onto
itself and satisfies conditions: $\psi(y_{1})=y_{1}$, and image of
$y_{0}$ is the point $y_{0}^{\perp}:=\psi(y_{0})$ opposite to $y_{0}$
on $S(x_{0},r)$. The coefficient of quasiconformality $Q_{1}$ of
$\psi$ can be easily estimated $Q_{1}\leq\frac{\pi r}{\bar{r}}=\frac{2}{\min(C_{r}-1,1)}$. 

By definition of quasiconformal homeomorphism we have 
\[
\begin{array}{c}
\mathrm{Cap}\left(F_{0}\cap B(x_{0},C_{r}r),F_{1}\cap B(x_{0},C_{r}r);B(x_{0},C_{r}r)\right)\geq\\
\frac{1}{Q_{1}}Cap\left(\psi(F_{0})\cap B(x_{0},C_{r}r),\psi(F_{1})\cap B(x_{0},C_{r}r);B(x_{0},C_{r}r)\right).
\end{array}
\]

Choose two closed balls 
$$B_{2}:=\bar{B}\left(x_{0}+\frac{y_{0}^{\perp}-x_{0}}{2},\frac{r}{2}\right)$$
and 
$$B_{3}:=\bar{B}\left(x_{0}+\frac{y_{0}^{\perp}-x_{0}}{2},\frac{r}{2}+(C_{r}-1)r\right)\subset B(x_{0},C_{r}r)$$
. By construction continua $\psi(F_{0})$ and $\psi(F_{1})$ intersect
any sphere $S\left(x_{0}+\frac{y_{0}^{\perp}-x_{0}}{2},\rho\right)$
for $\frac{r}{2}\leq\rho\leq\frac{r}{2}+(C_{r}-1)r$. Therefore by
Lemma 2 
\[
\mathrm{Cap}\left(\psi(F_{0})\cap B(x_{0},C_{r}r),\psi(F_{1})\cap B(x_{0},C_{r}r);B(x_{0},C_{r}r)\right)
\geq C(n)\ln\left(1+2(C_{r}-1)\right).
\]

Finally we have
\begin{equation}
\mathrm{Cap}\left(F_{0},F_{1};\Omega\right)\geq C(n)\min(C_{r}-1,1)\ln\left(1+2(C_{r}-1)\right).\label{eq:third}
\end{equation}

Combining this inequality and inequality \ref{eq:first} we obtain
finally
\begin{equation}
C(n)\min(C_{r}-1,1)\ln\left(1+2(C_{r}-1)\right)\leq\omega_{n-1}\left[\ln\left(K_{r}\right)\right]^{1-n}.\label{eq:final2}
\end{equation}

The proposition is proved for the second case:
 $$\mathrm{dist}(S_{0},S_{1})>\min\left(\left(C_{r}-1\right)\frac{r}{2},\frac{r}{2}\right)$$

Recall that $K_{r}$ is the short notation for $K_{int}(F_{r})$.
Combining inequalities (\ref{eq:first}, \ref{eq:second}, \ref{eq:third})
we obtain finally the constant $C_{1}(Q,C_{r},n)$ 
\[
C(n)\min\left[\ln2,\min(C_{r}-1,1)\ln\left(1+2(C_{r}-1)\right)\right]\leq\omega_{n-1}\left[\ln\left(K_{r}\right)\right]^{1-n}
\]
i.e.
\[
K_{r}\leq\exp\left[\left\{ \frac{\omega_{n-1}}{C(n)\min\left[\ln2,\min(C_{r}-1,1)\ln\left(1+2(C_{r}-1)\right)\right]}\right\} ^{\frac{1}{n-1}}\right].
\]

{\bf Remark.} {\it We proved a stronger estimate: 
\[
\frac{\bar{R}}{\bar{r}}\leq\exp\left[\left\{ \frac{\omega_{n-1}}{C(n)\min\left[\ln2,\min(C_{r}-1,1)\ln\left(1+2(C_{r}-1)\right)\right]}\right\} ^{\frac{1}{n-1}}\right].
\]}

Let us return to estimates of $\delta\left(F_{r}\right)$ with the help
of $\mathrm{dist}\left(F_{r},\partial\Omega'\right)$ i.e. to estimates
of $K_{ext}(F_{r})$. We start with the upper estimate:

{\bf Proposition 4.} {\it Suppose that $\Omega$,$\Omega'$ are domains
in $\mathbf{R}^{n}$ and the residual set $\Lambda=\mathbf{R}^{n}\setminus\Omega'$
of $\Omega'$ is unbounded and connected. Let $\varphi:\Omega\to\Omega'$
be a $Q$-quasiconformal homeomorphism of a bounded domain $\Omega\subset\mathrm{\mathbf{R}^{n}}$
onto a bounded domain $\Omega'\subset\mathrm{\mathbf{R}^{n}}$, $\bar{B}(x_{0},r)\subset\Omega$
be a closed ball with the embedding coefficient $K_{B(x_{0},r),\Omega}$.
 Then 
\[
\delta\left(F_{r}\right)\leq C_{1}\mathrm{(Q,C_{r},n)dist}\left(F_{r},\Lambda'\right)=C_{1}\mathrm{dist}\left(F_{r},\partial\Omega'\right)
\]
where a constant $C_{2}:=C_{2}(Q,C_{r},n)$ depends only on $Q,C_{r},n$.}

{\bf Proof.}
We will use notations: $F_{r}=\varphi(\bar{B}(x_{0},r))$$,K_{r}$
for the interior dilatation $K_{int}(F_{r})$, $C_{r}$ for the embedding
coefficient $K_{B(x_{0},r),\Omega}$ and $F_{R}:=\overline{\Omega}'\setminus\varphi(B(x_{0},C_{r}r))$. 

Because $\varphi:\Omega\to\Omega'$ is $Q$-quasiconformal
\[
\mathrm{Cap}(F_{r},F_{R};\Omega')\leq Q\:\mathrm{Cap}(\bar{B}(x_{0},r),\overline{\Omega}\setminus B(x_{0},C_{r}r));\Omega).
\]

By definition of the conformal capacity
\[
\:\mathrm{Cap}(\bar{B}(x_{0},r),\overline{\Omega}\setminus B(x_{0},C_{r}r));\Omega)=\:\mathrm{Cap}(\bar{B}(x_{0},r),\mathrm{\mathbf{R}^{n}}\setminus B(x_{0},Cr));\mathrm{\mathbf{R}^{n}}).
\]
By Lemma 1
\[
\mathrm{Cap}(\bar{B}(x_{0},r),\mathrm{\mathbf{R}^{n}}\setminus B(x_{0},Cr));\mathrm{\mathbf{R}^{n}})=\omega_{n-1}\left[\ln\left(C_{r}\right)\right]^{1-n}.
\]

Hence 
\[
\mathrm{Cap}(F_{r},F_{R};\Omega')\leq\omega_{n-1}\left[\ln\left(C_{r}\right)\right]^{1-n}
\]

Choose points $y_{o}\in F_{r}$and $y_{1}\in\partial\Omega'$ such
that $\mathrm{dist}(F_{r},\partial\Omega')=|y_{1}-y_{0}|$. Let $\tilde{y}:=y_{0}+\frac{y_{1}-y_{0}}{2}$
and $\bar{B(}\tilde{y},\rho_{1})$ be a closed ball with center at
$\tilde{y}$. 

By Proposition 3 $F_{r}$ contains a closed ball
of a radius $\rho_{2}:=\frac{1}{2}K_{r}\delta(F_{r})$. Therefore
both continuums $F_{r}$ and $F_{R}$ intersect any sphere $S(\tilde{y},\rho)$
for any $\rho\in[\rho_{1},\rho_{1}+\rho_{2}]$. By Lemma 2
we have
\[
\begin{array}{c}
\mathrm{Cap}(F_{r},F_{R};\Omega')\geq\mathrm{Cap}(F_{r},F_{R};\bar{B(}\tilde{y},\rho_{2})\setminus B(\tilde{y},\rho_{1}))\geq\\
C(n)\:\ln\left(\frac{\mathrm{dist}(F_{r},\partial\Omega')+K_{r}\delta(F_{r})}{\mathrm{dist}(F_{r},\partial\Omega')}\right).
\end{array}
\]
Combining both estimates for $\mathrm{Cap}(F_{r},F_{R};\Omega')$
we obtain finally 
\[
\mathrm{C(n)\:\ln\left(\frac{\mathrm{dist}(F_{r},\partial\Omega')+K_{r}\delta(F_{r})}{\mathrm{dist}(F_{r},\partial\Omega')}\right)\leq Cap}(F_{r},F_{R};\Omega')\leq\omega_{n-1}\left[\ln\left(C_{r}\right)\right]^{1-n}.
\]
It means that 
\[
\mathrm{dist}(F_{r},\partial\Omega')+K_{r}\delta(F_{r})\leq\mathrm{dist}(F_{r},\partial\Omega')\exp\left[\frac{\omega_{n-1}}{C(n)\left[\ln\left(C_{r}\right)\right]^{n-1}}\right].
\]
Because $C_{r}>1$ the number under exponent is positive and the exponent
value is bigger than $1$. Therefore 
\[
\delta(F_{r})\leq\left\{ K_{r}^{-1}\exp\left[\frac{\omega_{n-1}}{C(n)\left[\ln\left(C_{r}\right)\right]^{n-1}}\right]\right\} \mathrm{dist}(F_{r},\partial\Omega').
\]

The second estimate is an estimate from below with some additional restriction
on the domain $\Omega'$. We suppose additionally that any component of
$\mathrm{\mathbf{R}^{n}}\setminus\Omega'$ is unbounded. We call such
domains simple domains.

{\bf Proposition 5.} { \it Suppose that $\Omega$,$\Omega'$ are
domains in $\mathbf{R}^{n}$ and  the residual set $\Lambda=\mathbf{R}^{n}\setminus\Omega$
of $\Omega$ is unbounded and connected (i.e $\Omega$  is a simple domain). Let $\varphi:\Omega\to\Omega'$
be a $Q$-quasiconformal homeomorphism , $\bar{B}(x_{0},r)\subset\Omega$
be a closed ball with embedding coefficient $K_{B(x_{0},r),\Omega}$
. Then 
\[
\delta\left(F_{r}\right)\geq C_{3}\mathrm{(Q,K_{B(x_{0},r),\Omega},n)dist}\left(F_{r},\Lambda'\right)
\]
where a constant $C_{3}:=C_{3}(Q,K_{B(x_{0},r),\Omega},n)$ depends
only on \\ $Q,K_{B(x_{0},r),\Omega},n$.}

{\bf Proof.}
Denote by $\bar{B}_{r}:=\bar{B}(y_{0},\bar{r})\subset F_{r}$ a closed
ball of a maximal radius that belongs to $F_{r}:=\varphi(B(x_{o},r))$
and has the center at a point $y_{0};=\varphi(x_{0})$, by $B_{R}:=B(y_{0},\bar{R})$
an open ball of a minimal radius whose closure contains $F_{r}$,
by $CB^{R}:=\mathrm{\mathbf{R}^{n}}\setminus B^{R}$, by $K_{r}$
the ratio $\overline{K}_{r}\leq\frac{\bar{R}}{\bar{r}}$ and by $C_{r}$
the embedding coefficient $K_{B(x_{0},r),\Omega}$. By definition
of the interior dilatation $K_{r}\leq\overline{K}_{r}=\frac{\bar{R}}{\bar{r}}$. 

By Proposition 3 $\overline{K}_{r}\leq\tilde{C}_{r}(Q,C_{r},n):=\tilde{C}_{r}$.
It means that $\delta\left(F_{r}\right)\leq2\tilde{C}_{r}\overline{r}$.
Let $B_{r_{d}}:=B(y_{0},r_{d})$ be a greatest open ball such that
$S(y_{0},r_{d})\cap\partial\Omega'\neq\emptyset$ . By construction
$r_{d}\geq\overline{r}+dist\left(F_{r},\partial\Omega'\right)$. By
Lemma 1 and monotonicity of the conformal capacity
$$
\mathrm{Cap}(\bar{B}_{r},\Omega'\setminus B_{r_{d}};\mathrm{\Omega'})=\mathrm{Cap}(\bar{B}_{r},\mathrm{\mathbf{R}^{n}}\setminus B_{r_{d}};\mathrm{\mathbf{R}^{n}})=\omega_{n-1}\left[\ln\left(\frac{r_{d}}{\overline{r}}\right)\right]^{1-n}
$$
$$
\leq\omega_{n-1}\left[\ln\left(\frac{\overline{r}+dist\left(F_{r},\partial\Omega'\right)}{\overline{r}}\right)\right]^{1-n}.
$$

Using monotonicity of conformal capacity and quasi-invariance of conformal
capacity under $Q$-quasiconformal homeomorphisms we obtain
\[
\mathrm{\mathrm{Q^{-1}Cap}(\overline{B}(x_{0},r),\Lambda;\mathrm{\Omega})\leq\mathrm{Cap}(\bar{B}_{r},\Lambda';\mathrm{\Omega'})\leq\mathrm{Cap}(\bar{B}_{r},\mathrm{\mathbf{R}^{n}}\setminus B_{r_{d}};\mathrm{\Omega'}) \leq Cap}(\bar{B}_{r},\overline{\Omega}'\setminus B_{r_{d}};\mathrm{\Omega'}).
\]

Choose a point $x_{1}\in\partial\Omega$ closest to $\overline{B}(x_{0},r)$
and two balls 
$$B_{1}:=B\left(x_{0}+\frac{x_{1}-x_{0}}{2},\frac{C_{r}-1}{2}\right)$$
and 
$$\overline{B}_{2}:=\overline{B}\left(x_{0}+\frac{x_{1}-x_{0}}{2},\frac{C_{r}-1}{2}+2r\right)$$.
Continuums $F_{0}:=\overline{B}(x_{0},r)$ and $F_{1}:=\Lambda\cap B\left(x_{0}+\frac{x_{1}-x_{0}}{2},\frac{C_{r}-1}{2}r+2r\right)$
join spheres $S^{n-1}\left(x_{0}+\frac{x_{1}-x_{0}}{2},\frac{C_{r}-1}{2}r+2r\right)$
and $B\left(x_{0}+\frac{x_{1}-x_{0}}{2},\frac{C_{r}-1}{2}r\right)$.
By monotonicity of capacity and Lemma 2 we get
\[
\mathrm{C(n)\ln\left\{ 1+\frac{1}{C_{r}-1}\right\} \leq Cap}(F_{0},F_{1};\overline{B}_{2}\setminus B_{1})\leq\mathrm{Cap}(\overline{B}(x_{0},r),\Lambda;\mathrm{\Omega})
\]

Combining all previous estimates we get 
\[
\ln\left(1+\frac{1}{C_{r}-1}\right)\leq\frac{Q\omega_{n-1}}{C(n)}\omega_{n-1}\left[\ln\left(\frac{\overline{r}+dist\left(F_{r},\partial\Omega'\right)}{\overline{r}}\right)\right]^{1-n}.
\]

It means that 
\[
\left[\ln\left(1+\frac{dist\left(F_{r},\partial\Omega'\right)}{\overline{r}}\right)\right]^{n-1}\leq\frac{Q\omega_{n-1}}{C(n)}\omega_{n-1}\left[\ln\left(1+\frac{1}{C_{r}-1}\right)\right]^{-1}
\]

The constant 
\[
C_{0}(Q,C_{r},n):=\frac{Q\omega_{n-1}}{C(n)}\omega_{n-1}\left[\ln\left(1+\frac{1}{C_{r}-1}\right)\right]^{-1}
\]
is positive.

After elementary calculations we have 
\[
\frac{dist\left(F_{r},\partial\Omega'\right)}{\overline{r}}\leq\exp\left(C_{0}(Q,C_{r},n)^{\frac{1}{n-1}}\right)-1.
\]

Recall that $C_{r}$ is a short notation for $K_{B(x_{0},r),\Omega}$
and denote 
\[
C_{3}(Q,K_{B(x_{0},r),\Omega},n):=\left[\exp\left(C_{0}(Q,C_{r},n)^{\frac{1}{n-1}}\right)-1\right]^{-}.
\]

Because $\overline{r}\leq\delta(F_{r})$ we rewrite the last inequality
as 
\[
C_{3}(Q,K_{B(x_{0},r),\Omega},n)\:dist\left(F_{r},\partial\Omega'\right)\leq\overline{r}\leq\delta(F_{r})
\]

The main result of this paper can be formulated in more general way.

{\bf Theorem.}
{\it For any simple domain $\Omega$ quasiconformal image of any its Whitney family is a rough Whitney family.}

Follows directly from Propositions 3,4,5.

\address{Ben-Gurion University of the Negev, Department of Mathematics,\textbackslash newline
P.O. Box 653, 8410501 Beer Sheva, Israel.}
\email{vladimir@bgu.ac.il}

\address{College of William and Mary, Department of Mathematics,\textbackslash newline
P.O. Box 8795, Williamsburg, VA, USA.}
\email{nxzobi@wm.edu}

\end{document}